\overfullrule=0pt
\centerline {\bf More on the metric projection onto a closed convex set in a Hilbert space}\par
\bigskip
\bigskip
\centerline {BIAGIO RICCERI}\par
\bigskip
\bigskip
\bigskip
\bigskip
Here and in what follows, $(H,\langle\cdot,\cdot\rangle)$ is a real
Hilbert space and $X$ is a non-empty closed convex subset of $H$.
\smallskip
For each $x\in H$, we denote by $P(x)$ the metric projection of $x$ on $X$, that is the unique global minimum
of the restriction of the functional $y\to \|x-y\|$ to $X$.\par
\smallskip
There is no doubt that the map $P$ is among the most important and studied ones within convex analysis, functional analysis and 
optimization theory.\par
\smallskip
For the above reason, we think that it is of interest to highlight some properties of $P$ which do not appear in the wide literature concerning
$P$.\par
\smallskip
We collect such properties in Theorems 1, 2 and 3 below.\par
\smallskip
First, we fix some notations.\par
\smallskip
For each $r>0$, we put
$$B_r=\{x\in H : \|x\|^2<r\}$$
and
$$S_r=\{x\in H : \|x\|^2=r\}\ .$$
Moreover, for each $x\in X$, we set
$$J(x)={{1}\over {2}}(\|x\|^2-\|x-P(x)\|^2+\|P(0)\|^2)\ .$$
Furthermore, for each $r>0$, we put 
$$\gamma(r)=\inf_{x\in S_r}\|x-P(x)\|^2\ .$$
\medskip
Finally, since $P$ is non-expansive in $H$, for each $\lambda\in ]-1,1[$, the map $\lambda P$ is a contraction and hence has a unique fixed point that we 
denote by $\hat y_{\lambda}$.\par
\medskip
THEOREM 1. - {\it 
Assume that $0\not\in X$.\par
  Then, the following assertions hold:\par
\noindent
$(c_1)$\hskip 5pt the function $\lambda\to g(\lambda):=J(\hat y_{\lambda})$ is
increasing in $]-1,1[$ and its range is $\left ]-\|P(0)\|^2,\|P(0)\|^2\right [$\ ;\par
\noindent
$(c_2)$\hskip 5pt for each $r\in \left ]-\|P(0)\|^2,\|P(0)\|^2\right [$, the
point $\hat x_r:=\hat y_{g^{-1}(r)}$
is the unique point of minimal norm of $J^{-1}(r)$ towards which
every minimizing sequence in $J^{-1}(r)$, for the norm, converges\ ;  \par
\noindent
$(c_3)$\hskip 5pt the function $r\to \hat x_r$ is continuous in $\left ]-\|P(0)\|^2,\|P(0)\|^2\right [$\ ;\par
\noindent
$(c_4)$\hskip 5pt the function $\lambda\to h(\lambda):=\|\hat y_{{{1}\over {\lambda}}}\|^2$ is
decreasing in $]1,+\infty[$ and its range is $]0,\|P(0)\|^2[$\ ;\par
\noindent
$(c_5)$\hskip 5pt for each $r\in ]0, \|P(0)\|^2[$, the point $\hat v_r:=\hat y_{{{1}\over {h^{-1}(r)}}}$
is the unique global maximum of $J_{|S_r}$ towards which
 every maximizing sequence for $J_{|S_r}$ converges\ ;  \par
\noindent
$(c_6)$\hskip 5pt the function $r\to \hat v_r$ is continuous in $]0,\|P(0)\|^2[$\ .\par
\noindent
Assuming, in addition, that $X$ is compact, the following assertions hold:\par
\noindent
$(c_7)$\hskip 5pt the function $\gamma$ is $C^1$, decreasing and strictly convex in $]0,\|P(0)\|^2[$\ ; \par
\noindent
$(c_8)$\hskip 5pt
 one has
$$P(\hat v_r)=-\gamma'(r)\hat v_r$$
for all $r\in ]0,\|P(0)\|^2[$\ ;\par
\noindent
$(c_9)$\hskip 5pt one has
$$\gamma'(r)=-h^{-1}(r)$$
for all $r\in ]0,\|P(0)\|^2[$.}\par
\smallskip
PROOF. Clearly, the set of all fixed points of $P$ agrees with $X$. Now, fix $u\in H$ and $\lambda<1$. We show that
$$P(u+\lambda(P(u)-u))=P(u)\ .\eqno{(1)}$$
If $u\in X$, this is clear. Thus, assume $u\not\in X$ and hence $P(u)\neq u$.
Let $\varphi:H\to {\bf R}$ be the continuous linear functional defined by
$$\varphi(x)=\langle P(u)-u,x\rangle$$
for all $x\in H$. Clearly, $\|\varphi\|_{H^*}=\|P(u)-u\|$.
We have
$$\hbox {\rm dist}(u+\lambda(P(u)-u),\varphi^{-1}(\varphi(P(u)))=
{{|\varphi(u+\lambda(P(u)-u))-\varphi(P(u))|}\over {\|\varphi\|_{H^*}}}=(1-\lambda)\|P(u)-u\|\ .\eqno{(2)}$$
Moreover, by a classical result ([6], Corollary 25.23), we have
$$\langle P(u)-u,P(u)-x\rangle\leq 0$$
for all $x\in X$, that is
$$X\subseteq \varphi^{-1}([\varphi(P(u)),+\infty[)\ .\eqno{(3)}$$
Also, notice that 
$$\hbox {\rm dist}(u+\lambda(P(u)-u),\varphi^{-1}(\varphi(P(u))))=
\hbox {\rm dist}(u+\lambda(P(u)-u),\varphi^{-1}([\varphi(P(u)),+\infty[))\ . \eqno{(4)}$$
Indeed, otherwise, it would exist $w\in H$, with $\varphi(w)>\varphi(P(u))$, such that
$$\|u+\lambda(P(u)-u)-w\|<\hbox {\rm dist}(u+\lambda(P(u)-u),\varphi^{-1}(\varphi(P(u)))\ .$$
Then, since $\varphi(u+\lambda(P(u)-u))<\varphi(P(u))$ (indeed $\varphi(u+\lambda(P(u)-u))-\varphi(P(u))=(\lambda-1)\|P(u)-u\|^2$),
by connectedness and continuity, in the open ball centered at $u+\lambda(P(u)-u)$, of radius 
$\hbox {\rm dist}(u+\lambda(P(u)-u),\varphi^{-1}(\varphi(P(u)))$, it would exists a point at which $\varphi$ takes the value
$\varphi(P(u))$, which is absurd. So, $(4)$ holds. Now, from $(2), (3), (4)$, it follows that
$$(1-\lambda)\|P(u)-u\|\leq \hbox {\rm dist}(u+\lambda(P(u)-u),X)\leq \|u+\lambda(P(u)-u)-P(u)\|=(1-\lambda)\|P(u)-u\|$$
which yields $(1)$. From $(1)$, in particular, we infer that $P(0)=P(-P(0))$. On the other hand, if $\tilde x\in H$
is such that $\tilde x=-P(\tilde x)$, then,
applying $(1)$ with $u=\tilde x$ and $\lambda={{1}\over {2}}$, we get $P(0)=P(\tilde x)$
and so $\tilde x=-P(0)$. Therefore, $-P(0)$ is the unique fixed point of $-P$. Now, let us recall that
$J$ is a Fr\'echet differentiable convex functional whose derivative is equal to $P$ ([1], Proposition 2.2).
This allows us to use the results of [3].  Therefore,  $(c_1)$, $(c_2)$, $(c_3)$ follow respectively from $(a_1)$, $(a_2)$, $(a_3)$ of Theorem
3.2 of [3], since (with the notation of that result) we have $\eta_1=J(-P(0))=-\|P(0)\|^2$ and $\theta_1=\inf_XJ=\|P(0)\|^2$, while $(c_4)$, $(c_5)$,
$(c_6)$ follow respectively from $(b_1)$, $(b_2)$, $(b_3)$ of Theorem 3.3 of [3], since $\theta_2=\|P(0)\|^2$.
Now, assume that
$X$ is also compact. Then, $J$ turns out to be sequentially weakly continuous ([5], Corollary 41.9). Moreover, $J$ has no local maxima since $P$ 
has no zeros.
At this point, $(c_7)$, $(c_8)$, $(c_9)$  follow respectively from $(b_4)$, $(b_5)$, $(b_6)$ of Theorem 3.3 of [3], since, for a constant $k_0$,
we have 
$$\sup_{S_r}J=-{{1}\over {2}}\gamma(r) +k_0 $$
for all $r>0$. The proof is complete.\hfill $\bigtriangleup$
\medskip
THEOREM 2. - {\it Let $Q:H\to H$ be a continuous and monotone potential operator such that
$$\lim_{\|x\|\to +\infty}I(x):=\int_0^1\langle Q(sx),x\rangle ds=+\infty\ .$$
Set
$$\lambda^*=\inf_{r>\inf_HI}\inf_{x\in I^{-1}(]-\infty,r[)}{{J(x)-\inf_{y\in I^{-1}(]-\infty,r])}J(y)}\over {r-I(x)}}\ .$$
Then,  the equation
$$P(x)+\lambda Q(x)=0$$
has a solution in $H$ for every $\lambda>\lambda^*$. Moreover, when $\lambda^*>0$, the same equation has no solution in $H$ 
for every $\lambda<\lambda^*$.}\par
\smallskip
PROOF. Since $Q$ is a monotone potential operator, the functional $I$ turns out be convex, of class $C^1$ and its derivative agrees with $Q$. Now, the 
conclusion
follows from Theorem 2.4 of [2], since, by convexity, the solutions of the equation $P(x)+\lambda Q(x)$ are exactly the global minima in $H$ of the 
functional $J+\lambda I$.
\hfill $\bigtriangleup$
\medskip
THEOREM 3. - {\it Let $(T, {\cal F},\mu)$ be a  measure space, with $0<\mu(T)<+\infty$ and assume that $0\not\in X$.\par
Then, for every $\eta\in L^{\infty}(T)$, with $\eta\geq 0$, for every $r\in ]0,\|P(0)\|^2[$ and for every $p\geq 2$, if we put
$$U_{\eta,r}=\left \{u\in L^p(T,H) :\int_T\eta(t)\|u(t)\|^2d\mu= r\int_T\eta(t)d\mu\right\}\ ,$$
we have
$$\inf_{u\in U_{\eta,r}}\int_T\eta(t)\|u(t)-P(u(t))\|^2d\mu
=\inf_{x\in S_r}\|x-P(x)\|^2\int_T\eta(t)d\mu \eqno{(5)}$$
and
$$\sup_{u\in U_{\eta,r}}\int_T\eta(t)\|u(t)-P(u(t))\|^2d\mu
=\sup_{x\in S_r}\|x-P(x)\|^2\int_T\eta(t)d\mu \ .\eqno{(6)}$$}
\smallskip
PROOF. Applying Theorem 5 of [4] to $J$ and $-J$, respectively, we obtain
$$\inf_{u\in V_{\eta,r}}\int_T\eta(t)J(u(t))d\mu=\inf_{S_r}J\int_T\eta(t)d\mu \eqno{(7)}$$
and
$$\sup_{u\in V_{\eta,r}}\int_T\eta(t)J(u(t))d\mu=\sup_{ S_r}J\int_T\eta(t)d\mu \ ,\eqno{(8)}$$
where
$$V_{\eta,r}=\left \{u\in L^p(T,H) :\int_T\eta(t)\|u(t)\|^2d\mu\leq  r\int_T\eta(t)d\mu\right\}\ .$$
Now, observe that $J_{|S_r}$ has a global minimum. Indeed, since $J$ is weakly lower semicontinuous and $\overline {B_r}$ is weakly
compact, $J_{|\overline {B_r}}$ has a global minimum, say $\hat w_r$. Notice that $\hat w_r\in S_r$, since, otherwise, $P(\hat w_r)=0$
which is impossible since $0\not\in X$. So, $\hat w_r$ is a global minimum of $J_{|S_r}$. Furthermore, from Theorem 1, we know that
$J_{|S_r}$ has a global maximum, say $\hat v_r$. Denote by the same symbols the constant functions (from $T$ into $Y$) taking, respectively, 
the values $\hat w_r$ and $\hat v_r$. Since $\mu(T)<+\infty$, we have $\hat w_r, \hat v_r\in U_{\eta,r}$. So, from $(7)$ and $(8)$, it follows
respectively
$$\inf_{u\in V_{\eta,r}}\int_T\eta(t)J(u(t))d\mu=\int_TJ(\hat w_r)\eta(t)d\mu\geq \inf_{u\in U_{\eta,r}}\int_T\eta(t)J(u(t))d\mu$$
and 
$$\sup_{u\in V_{\eta,r}}\int_T\eta(t)J(u(t))d\mu=\int_TJ(\hat v_r)\eta(t)d\mu\leq \sup_{u\in U_{\eta,r}}\int_T\eta(t)J(u(t))d\mu\ .$$
Therefore
$$\inf_{S_r}J\int_T\eta(t)d\mu=(r+\|P(0)\|^2-\sup_{x\in S_r}\|x-P(x)\|^2)\int_T\eta(t)d\mu=\inf_{u\in U_{\eta,r}}
\int_T\eta(t)(\|u(t)\|^2-\|u(t)-P(u(t))\|^2+\|P(0)\|^2)d\mu$$
$$=(r+\|P(0)\|^2)\int_T\eta(t)d\mu-\sup_{u\in U_{\eta,r}}\int_T\eta(t)
\|u(t)-P(u(t))\|^2d\mu$$
which yields $(6)$. Likewise
$$\sup_{S_r}J\int_T\eta(t)d\mu=(r+\|P(0)\|^2-\inf_{x\in S_r}\|x-P(x)\|^2)\int_T\eta(t)d\mu=\sup_{u\in U_{\eta,r}}
\int_T\eta(t)(\|u(t)\|^2-\|u(t)-P(u(t))\|^2+\|P(0)\|^2)d\mu$$
$$=(r+\|P(0)\|^2)\int_T\eta(t)d\mu-\inf_{u\in U_{\eta,r}}\int_T\eta(t)
\|u(t)-P(u(t))\|^2d\mu$$
which yields $(5)$\hfill $\bigtriangleup$
\vfill\eject
\centerline {\bf References}\par
\bigskip
\bigskip
\noindent
[1]\hskip 5pt S. FITZPATRICK and R. PHELPS, {\it Differentiability of the metric projection in Hilbert space},
 Trans. Amer. Math. Soc., {\bf 270} (1982), 483-501.\par
\smallskip
\noindent
[2]\hskip 5pt B. RICCERI, {\it A general variational principle and
some of its applications}, J. Comput. Appl. Math., {\bf 113}
(2000), 401-410.\par
\smallskip
\noindent
[3]\hskip 5pt B. RICCERI, {\it Fixed points of nonexpansive potential
operators in Hilbert spaces}, Fixed Point Theory Appl. {\bf 2012},  2012: 123.\par
\smallskip
\noindent
[4]\hskip 5pt B. RICCERI, {\it Integral functionals on $L^p$-spaces:
infima over sub-level sets}, Numer. Funct. Anal. Optim., {\bf 35} (2014), 1197-1211.\par
\smallskip
\noindent
[5]\hskip 5pt E. ZEIDLER, {\it Nonlinear Functional Analysis and its Applications}, vol. III, Springer-Verlag, 1985.\par
\smallskip
\noindent
[6]\hskip 5pt E. ZEIDLER, {\it Nonlinear Functional Analysis and its Applications}, vol. II/B, Springer-Verlag, 1990.\par
\bigskip
\bigskip
\bigskip
\bigskip
Department of Mathematics\par
University of Catania\par
Viale A. Doria 6\par
95125 Catania, Italy\par
{\it e-mail address}: ricceri@dmi.unict.it

\bye